\documentclass[11pt,a4]{article}
\usepackage{amsmath}
\usepackage{amssymb}
\usepackage{theorem}
\usepackage{color}
\newcommand{\R}{{\cal R}}
\newcommand{\B}{{\cal B}}
\newcommand{\F}{{\cal F}}

\newcommand{\Hc}{{\cal H}}
\newcommand{\Nb}{{\mathbb{N}}}
\newcommand{\Qb}{{\mathbb{Q}}}
\newcommand{\Rb}{{\mathbb{R}}}
\newcommand{\Zb}{{\mathbb{Z}}}

\newcommand{\qed}{\hfill\hbox{\rule{6pt}{6pt}}}
\newtheorem{theorem}{Theorem}[section]
\newtheorem{lemma}{Lemma}[section]
\newtheorem{col}{Corollary}[section]

{\theorembodyfont{\rmfamily}
{\theorembodyfont{\rmfamily}

\begin{document}
\title{Generalization of van Lambalgen's  theorem and blind randomness for conditional probabilities}
\author{Hayato Takahashi\thanks{Tokyo Denki University, Department of  Electrical and  Electronic Engineering,  5 Senju Asahi-cho, Adachi-ku, Tokyo 120-8551 Japan.
Present address: Gifu University, Japan.   hayato.takahashi@ieee.org}}
\date{\today}
\maketitle

\begin{abstract}
Generalization of the Lambalgen's theorem  is studied with
 the notion of Hippocratic (blind) randomness without assuming computability of conditional probabilities. 
In [Bauwence 2014], a counter-example for the generalization of Lambalgen's theorem is shown when the conditional probability is not computable.
In this paper, it is shown that (i) finiteness of martingale for blind randomness, (ii) classification of two blind randomness by likelihood ratio test, (iii) sufficient conditions for 
the generalization of the Lambalgen's theorem, and (iv) an example that satisfies the Lambalgen's theorem but the conditional probabilities are not computable for all random parameters. 
\end{abstract}
\section{Introduction}
Lambalgen's theorem (1987) \cite{lambalgen87} says that a pair of sequences \((x^\infty,y^\infty)\in\Omega^2\) is Martin-L\"of (ML) random w.r.t.~the product measure of uniform measures iff \(x^\infty\) is ML-random and \(y^\infty\) is ML-random relative to \(x^\infty\),
where \(\Omega\) is the set of infinite binary sequences. 
In this paper we study a generalized form of the Lambalgen's theorem using the notion of blind (Hippocratic) randomness \cite{{Laurent-etal},{kjoshanssen}}.

Let \(S\) be the set of finite binary strings and \(\Delta(s):=\{sx^\infty \vert x^\infty\in\Omega\}\) for \(s\in S\), where \(sx^\infty\) is the concatenation of \(s\) and \(x^\infty\).
Let \(|s|\) be the length of \(s\in S\). 
Let \(\Delta(x,y):=\Delta(x)\times\Delta(y)\) for \(x,y\in S\).
In this paper, we study probabilities on \((\Omega, \B)\) or \((\Omega^2, \B_{\Omega^2})\), where \(\B\) and \(\B_{\Omega^2}\) are the  Borel-\(\sigma\)-algebras generated from \(\Delta(s), s\in S\) and \(\Delta(x,y), (x,y)\in S^2\), respectively. 
In the following we omit \(\B\) or \(\B_{\Omega^2}\) and  write such as  \(P\) on \(\Omega\) or \(\Omega^2\), when it is obvious from the context. 
For a probability \(P\) on \(\Omega\), we write \(P(s):=P(\Delta(s))\) for \(s\in S\).
For a probability \(P\) on \(X\times Y, X=Y=\Omega\), let \(P(x,y):=P(\Delta(x,y))\) for \(x,y\in S\). Let  \(P_X, P_Y\) be the marginal distributions on \(X\) and \(Y\), respectively, i.e.,
\(\forall x,y\ \ P_X(x):=P(x,\lambda)\text{ and }P_Y(y):=P(\lambda,y)\).
For \(A\subseteq S\) and \(B\subseteq S^2\), set \(\tilde{A}:=\cup_{s\in A}\Delta(s)\) and \(\tilde{B}:=\cup_{(x,y)\in B}\Delta(x,y)\), respectively. 
For \(x,y\in S\), we write \(x\sqsubseteq y\) if \(x\) is a prefix of \(y\). Let \(\Nb\) be the set of natural numbers. 
Let \(\R^P\) be the set of ML-random set w.r.t.~\(P\) when \(P\) is computable and  \(\R^{P(\cdot\vert y^\infty), y^\infty}\) 
be the set of ML-random set w.r.t.~\(P(\cdot\vert y^\infty)\) relative to \(y^\infty\) when \(P(\cdot\vert y^\infty)\) is computable relative to \(y^\infty\), respectively.

For \(A\subseteq \Omega^2\text{ and } y^\infty\in\Omega\), set   \(A_{y^\infty}:=\{x^\infty\mid (x^\infty,y^\infty)\in A\}\).
For example, if \(P\) is a computable probability on \(\Omega^2\),  we write \(\R^P_{y^\infty}:=\{ x^\infty\mid (x^\infty, y^\infty)\in\R^P\}\) for \(y^\infty\in\Omega\).

In Vovk and Vyugin (1993) \cite{vovkandvyugin93}, they generalized Lambalgen's theorem as follows (actually they show a different form of the  following theorem with parametric models, however the following form is easily derived from them).
\begin{theorem}[Vovk and Vyugin \cite{vovkandvyugin93}]\label{general-vovk}
Let \(P\) be a computable probability on \(X\times Y, X=Y=\Omega\).
Assume that\\
 (i) conditional probabilities exist for all parameters and\\
  (ii) they are uniformly computable
for all parameters. \\
Then 
\[ \R^P_{y^\infty}=\R^{P(\cdot\vert y^\infty),y^\infty}\text{ for all }y^\infty\in\R^{P_Y}.\]
\end{theorem}

Here conditional probability \(P(\cdot\vert\cdot)\)   is called uniformly computable for all parameters if\\
 (i) there is a partial computable function \(A\) such that \\
\( \forall s\in S, y^\infty\in\Omega, k\in\Nb\exists y\sqsubset y^\infty,\ \vert P(s\vert y^\infty)-A(s,y,k)\vert <\frac{1}{k}\) 
and  \\
(ii) if \(A(s,y,k)\) is defined then \(A(s,y,k)=A(s, y',k)\) for all \(y'\sqsupseteq y\).\\
It is known that  there are non-uniform computable conditional probabilities (Roy 2011 \cite{roy2011}). 

In \cite{takahashiIandC}, it is shown that conditional probabilities exist for all random parameters, i.e., 
\begin{equation}\label{def-cond-prob}
\forall x\in S,\ y^\infty\in\R^{P_Y}\ \ P(x\vert y^\infty):=\lim_{y\to y^\infty}P(x\vert y)\text{ (the right-hand-side exist)}
\end{equation}
and \(P(\cdot\vert y^\infty)\) is a probability on \((\Omega, \B)\) for each \(y^\infty\in\R^{P_Y}\).

For any fixed \(y^\infty\in\Omega\),  \(P(\cdot\vert y^\infty)\) is called computable relative to \(y^\infty\) if \\
(i) there is a partial computable function \(A\) such that \\
\( \forall s\in S, k\in\Nb\exists y\sqsubset y^\infty,\ \vert P(s\vert y^\infty)-A(s,y,k)\vert <\frac{1}{k}\) 
and \\
 (ii) if \(A(s,y,k)\) is defined then \(A(s,y,k)=A(s, y',k)\) for all \(y'\sqsupseteq y\). 
 
 Note that in this definition, \(A\) may depend on \(y^\infty\), however in the definition of 
uniform computability of conditional probability,  we require that there is a global \(A\) that satisfies the conditions for all \(y^\infty\). 

The next theorem shows that  the generalized Lambalgen's theorem holds if the conditional probability is computable relative to the given parameter.
\begin{theorem}[\cite{takahashiIandC, takahashiIandC2}]\label{takahashi}
Let \(P\) be a computable probability on \(X\times Y, X=Y=\Omega\).
Fix \(y^\infty\in\R^{P_Y}\) and assume that the conditional probability \(P(\cdot\vert y^\infty)\) is computable relative to \(y^\infty\).
Then
\begin{equation}\label{T}
\R^P_{y^\infty}=\R^{P(\cdot\vert y^\infty),y^\infty}.
\end{equation}
\end{theorem}

Conditional probabilities always exist for all random parameters, however they may not be computable, see Theorem~\ref{exampleA} and \ref{exampleB} below. 
In this paper, we introduce the notion of blind (Hippocratic)  randomness and study the  generalization of (\ref{T}) when the conditional probability is not computable. 
Here blind randomness \(\Hc\) is defined as follows. 
Let \(P\) be a probability on \(\Omega\). 
An r.e.~set \(U\subseteq\Nb\times S\) is called blind test w.r.t.~\(P\) if \(U_n\supseteq U_{n+1}\) and \(P(\tilde{U}_n)<2^{-n}\), where \(U_n:=\{ x\mid (n,x)\in U\}\), for all \(n\).
The set of blind random sequences w.r.t.~\(P\) (in the following we denote it by \(\Hc^P\)) is the set that is not covered by any limit of blind test, i.e.,  \(\Hc^P:=(\cup_{U: \text{blind test}}\cap_n \tilde{U}_n)^c\) \cite{{Laurent-etal},{kjoshanssen}}.
Similarly a blind test \(U^{y^\infty}\) w.r.t.~\(P(\cdot\vert y^\infty)\) relative to \(y^\infty\)  is an r.e.~set relative to \(y^\infty\) such that \(U_n^{y^\infty}\supseteq U_{n+1}^{y^\infty}\) and \(P(\tilde{U}_n^{y^\infty}\vert y^\infty)<2^{-n}\), 
where \(U_n^{y^\infty}:=\{ x\mid (n,x)\in U^{y^\infty}\}\), for all \(n\).  Let \(\Hc^{P(\cdot\vert y^\infty),y^\infty}\) be the set of blind random sequences w.r.t.~\(P(\cdot\vert y^\infty)\) relative to \(y^\infty\), i.e., 
\(\Hc^{P(\cdot\vert y^\infty),y^\infty}\) is the set that is not covered by any limit of blind test w.r.t.~the conditional probability
relative to \(y^\infty\).  
If a probability is not computable, the existence of the universal test is not assured, however the  definitions  above are still well defined.
If \(P\) is computable, we have \(\R^P=\Hc^P\), and if \(P(\cdot\vert y^\infty)\) is computable relative to \(y^\infty\), we have \(\R^{P(\cdot\vert y^\infty),y^\infty}=\Hc^{P(\cdot\vert y^\infty),y^\infty}\).
In the definition above, we can replace \(P(\tilde{U}_n)<2^{-n}\) and \(P(\tilde{U}_n^{y^\infty}\vert y^\infty)<2^{-n}\) with \(P(\tilde{U}_n)<f(n)\) and \(P(\tilde{U}_n^{y^\infty}\vert y^\infty)<f(n)\), respectively, where \(f\) is a computable decreasing function.

In \cite{takahashiIandC}, in the proof of \(\supseteq\) part in Theorem~\ref{takahashi}, computability of conditional probability is not assumed. 
\begin{col}[\cite{takahashiIandC}]\label{col-ifpart}
Let \(P\) be a computable probability on \(\Omega^2\). Then 
\begin{equation}\label{generalB}
\R^P_{y^\infty}\supseteq\Hc^{P(\cdot\vert y^\infty),y^\infty} \text{ for all }y^\infty\in\R^{P_Y}.
\end{equation}
\end{col}
\section{Results}
First we show a sufficient condition for the equality in (\ref{generalB}).
In the following we set  \(\frac{a}{0}:=\infty\) if \(a\ne 0\) else \(0\).
\begin{theorem}\label{th-main}
Let \(P\) be a computable probability on \(\Omega^2\). 
Fix a pair of sequences \((x^\infty,y^\infty)\in\Omega^2\).
Assume that   there are a computable probability \(Q\) on \(\Omega^2\)  and a partial computable function with oracle \(y^\infty\), \(f_{y^\infty}: S\to \Qb\) 
such that
\begin{enumerate}
\item   \(y^\infty\in\R^{P_Y}\cap\R^{Q_Y}\),
\item  \(Q(\cdot\vert y^\infty)\) is computable relative to \(y^\infty\), 
\item \(\forall x\sqsubset x^\infty\ P(x\vert y^\infty)>0\), 
\item There is an infinite subset \( A\subseteq \{ x\mid x\sqsubset x^\infty\}\) such that \\
\(\sup_{x\in A}f_{y^\infty}(x)<\infty\) and  \\
\(A\subseteq \{ x\mid  f_{y^\infty}(x)\text{ is defined }\}\subseteq \{ x\mid  \frac{Q(x\vert y^\infty)}{P(x\vert y^\infty)}< f_{y^\infty}(x)<\infty\}\), and
\item
\(0<\inf_{x\to x^\infty}\frac{Q(x\vert y^\infty)}{P(x\vert y^\infty)}\).
\end{enumerate}
Then 
\[x^\infty\in \R^P_{y^\infty}\Leftrightarrow x^\infty\in \Hc^{P(\cdot\vert y^\infty),y^\infty}.\]
\end{theorem}
Proof)
The proof is almost same with  the proof of Theorem~3.3 in \cite{takahashiIandC2}.
Fix \((x^\infty,y^\infty)\) that satisfies the condition of the theorem. 
As in the proof of Theorem~3.3 in \cite{takahashiIandC2}, we expand a test w.r.t.~\(P(\cdot\vert y^\infty)\) to a global test w.r.t.~\(P\).
The problem here is that we do not assume the computability of \(P(\cdot\vert y^\infty)\).
However from the condition of the theorem, we can approximate the conditional probability with some computable function as follows.

From (iv) and (v),  let \(c_1\) and \(c_2\) be  rational constants such that 
\begin{equation}\label{main-bound}
0<c_1<\inf_{x\to x^\infty}\frac{Q(x\vert y^\infty)}{P(x\vert y^\infty)}\text{ and }\sup_{x\in A}f_{y^\infty}(x)<c_2<\infty.
\end{equation}

Let \(U\subseteq S\) be an r.e.~set relative to \(y^\infty\) such that 
\[
P(\tilde{U}\vert y^\infty)<2^{-n}c^{-1}_2\text{ and }x^\infty\in \tilde{U}.
\]
Let
\[V:=\{x\mid \exists z\in U\ z\sqsubseteq x,\ f_{y^\infty}(x)<c_2\}\text{ and }V':=\{x\in U\mid \frac{Q(x\vert y^\infty)}{P(x\vert y^\infty)}<c_2\}.\]
From (iv), we have \(x^\infty\in\tilde{V}\),  \(V\) is r.e.~relative to \(y^\infty\), and \(\tilde{V}\subseteq\tilde{V}'\).
From (iv) and (\ref{main-bound}), we have
\[Q(\tilde{V}\vert y^\infty)\leq Q(\tilde{V}'\vert y^\infty)<c_2P(\tilde{U}\vert y^\infty)<2^{-n}.\]

From Theorem~3.3 in \cite{takahashiIandC2}
\footnote{there is a typos error in the proof of Theorem~3.3 in \cite{takahashiIandC2}; (i) \(k\geq 2^{n+|x|}\) should be changed to \(k\geq 2^{n+|x|+1}\) in equation (7) pp.188,  line 1 and 3  pp.189., and 
(ii) inequality \(<\) should be changed to \(\leq\) in equation (15) pp188,   line 1 and 7 from bottom pp188.}, 
there is an r.e.~set \(W\subseteq S\times S\) such that 
\(Q(\tilde{W})<\frac{11}{2}2^{-n}\) and \(\tilde{W}_{y^\infty}=\tilde{V}\).
Let \(W':=\{(x,y)\in W\mid P(x\vert y)<Q(x\vert y)c_1^{-1}\}\).
Then \(W'\) is an r.e.~set, and from (\ref{def-cond-prob}) and  (\ref{main-bound}), we have
 \[
 P(\tilde{W'})<Q(\tilde{W'})c_1^{-1}<\frac{11}{2}2^{-n}c_1^{-1},\text{ and }x^\infty\in\tilde{W'}_{y^\infty}.
 \]

Therefore if \(x^\infty\) is covered by a test w.r.t.~\(P(\cdot\vert y^\infty)\) then \((x^\infty,y^\infty)\) is covered by a test w.r.t.~\(P\), which shows only if part of the theorem. 
The if part follows from Corollary~\ref{col-ifpart}.
\qed

\subsection{Martingale and likelihood ratio test for blind randomness}
Next we show a classification of  blind randomness for two different probabilities.
For similar results for ML-randomness, see  \cite{{BM2007},{takahashiIandC2}}.
Let \(P\) and \(Q\) be probabilities on \(\Omega\).
From martingale convergence theorem, we have
\[\lim_{x\to x^\infty} \frac{Q(x)}{P(x)}<\infty,\ a.s.-P.\]
In \cite{takahashiIandC2}, it is shown that martingale convergence theorem holds for individual ML-random sequences for computable \(P\), and the above inequality holds for them. 
In order to explore similar results for blind randomness without assuming computability of probabilities, we introduce a notion of approximation. 

Let \(\Zb\) be the set of the integers and \(\Rb\) be the set of real numbers, respectively. 
Let \(\F_n\) be the algebra generated from \(\{\Delta(x)\mid |x|=n\}\) then we have \(\B=\sigma(\cup_n\F_n)\).
Let \(r_n:\Omega\to\Rb\) be a measurable function w.r.t.~\(\F_n\). Since \(r_n(x^\infty)\) takes a constant value on  \(\Delta(x)\) for \( |x|=n, x\sqsubset x^\infty\),
we write 
\begin{equation}\label{def-r}
\forall n\in\Nb, x\in S, x^\infty\in\Omega\ \  r_n(x):=r_n(x^\infty)\text{ if }|x|=n, x\sqsubset x^\infty.
\end{equation}
Let \(g:\Rb \to\Rb\) be a strictly increasing function, i.e., \(\forall x,y\in\Rb\ g(x)<g(y)\Leftrightarrow x<y\).
We say that the set of random variables \(\{r_n\}_{n\in\Nb}\) is \(g\){\it-effectively-approximable} if there is a computable \(f:S\times\Nb\to\Qb\) such that 
\begin{equation}\label{r-approx}
\exists c\in\Nb\forall n\in\Nb\forall x\in S\ f(x,n)\leq g(r_n(x))\leq f(x,n)+c.
\end{equation}
We say that
\begin{enumerate}
\item
\(\{r_n\}_{n\in\Nb}\) is  {\it effectively-approximable} if there is a strictly increasing \(g\) such that \(\{r_n\}_{n\in\Nb}\) is \(g\)-effectively-approximable, 
\item
\(\{r_n\}_{n\in\Nb}\) is  {\it strongly-effectively-approximable} if  there is a strictly increasing \(g\) such that \(\{r_n\}_{n\in\Nb}\) is \(g\)-effectively-approximable 
and the restriction \(g\) to \(\Nb\) is computable, i.e.,  if \(\forall n\in\Nb\ g'(n):=g(n)\) then \(g'\) is computable,
\item
 \(\{r_n\}_{n\in\Nb}\) is  submartingale (w.r.t.~\(P\)) if \(\forall n\ E(r_n\vert \F_{n-1})\geq r_{n-1}\  P-a.s.\),  
 \item
\(\{r_n\}_{n\in\Nb}\) is  martingale (w.r.t.~\(P\)) if \(\forall n\ E(r_n\vert \F_{n-1})= r_{n-1}\  P-a.s.\).
\end{enumerate}

\begin{theorem}\label{martingale-blind}
Let \(P\) be a probability on \(\Omega\).
Let  \(\{r_n\}_{n\in\Nb}\) be a  non-negative submartingale w.r.t.~\(P\).
If \(\sup_i E(|r_i|)<\infty\) and \(\{r_n\}_{n\in\Nb}\) is strongly-effectively-approximable, then 
\[\Hc^P\subseteq \{x^\infty\mid \sup_n r_n(x^\infty)<\infty\}.\]
\end{theorem}
Proof)
Let 
\begin{align*}
U_{m,n} &:=\{x\mid m<\sup_{1\leq i\leq n}g(r_i(x))\}, \\
U'_{m,n} &:=\{x\mid m<g(\sup_{1\leq i\leq n}r_i(x))\}, \\
V_{m,n}&:=\{x\mid m<\sup_{1\leq i\leq n}f(x,i)\}, \text{ and}\\
M_{m,n}&:= \{x\mid m<\sup_{1\leq i\leq n}r_i(x)\}\text{ for all }m,n\in\Nb.
 \end{align*}
 Let 
\[V_{g(m)}:=\cup_n V_{g(m),n}=\{x\mid g(m)<\sup_{i\in\Nb} f(x,i)\}.\]
Since \(f\) and \(g\) are computable, we see that \(\{V_{g(m)}\}_{m\in\Nb}\) is an uniformly r.e.~set.
Since \(g\) is increasing, from (\ref{r-approx}), we have
\begin{equation}\label{marting-eq-A}
\exists c\forall x,n\  \sup_{1\leq i\leq n} f(x,i)\leq \sup_{1\leq i\leq n}g(r_i(x))=g(\sup_{1\leq i\leq n}r_i(x))\leq\sup_{1\leq i\leq n} f(x,i)+c, \text{ and }
\end{equation}
\begin{equation}\label{marting-eq-B}
\exists c\forall x\  \sup_{i\in\Nb} f(x,i)\leq g(\sup_{i\in\Nb}r_i(x))\leq\sup_{i\in\Nb} f(x,i)+c.
\end{equation}

Then 
\begin{align}
P(\tilde{V}_{g(m),n})& \leq P(\tilde{U}_{g(m),n})\label{martingale-1}\\
& = P(\tilde{U}^{\prime}_{g(m),n})\label{martingale-2}\\
& = P(\tilde{M}_{m,n})\label{martingale-3}\\
&\leq E(|r_n|)/m\label{martingale-4},
\end{align}
where (\ref{martingale-1}) follows from (\ref{marting-eq-A}), (\ref{martingale-2}) and (\ref{martingale-3}) follows from that \(g\) is strictly increasing, 
and  (\ref{martingale-4}) follows from Doob's submartingale inequality (for example, see \cite{williams91} pp.137).

Thus we have
\[P(\tilde{V}_{g(m)})=\lim_n P(\tilde{V}_{g(m),n})\leq \sup_n E(|r_n|)/m.\]
Since \(\sup_n E(|r_n|)<\infty\), we see that \(\{V_{g(m)}\}_{m\in\Nb}\) is a test. 
Let 
\[V'_m:=\{x\mid m<\sup_{i\in\Nb} r_i(x)\}.\]
Since \(g\) is strictly increasing, from (\ref{marting-eq-B}), we have
\(\cap_m \tilde{V}_{g(m)}=\cap_m \tilde{V}'_m\subseteq (\Hc^P)^c\).
\qed

For example, \(\frac{Q}{P}\) is {\it \(\log\)-effective-approximable} if there is a computable function  \(f:S\to \Zb\)  such that 
\(\exists c \forall x\ f(x)<\log\frac{Q(x)}{P(x)}<f(x)+c\).

\begin{col}\label{col-q/p}Let \(P\) and \(Q\) be probabilities on \(\Omega\).
 Suppose that \(\frac{Q}{P}\) is strongly-effective-approximable. Then
   \[\Hc^P\subseteq \{ x^\infty\mid \sup_{x\to x^\infty}\frac{Q(x)}{P(x)}<\infty\}.\]
\end{col}
Proof)
Since \(E(Q/P)\leq 1\), where the expectation is taken w.r.t.~\(P\), from Theorem~\ref{martingale-blind}, we have the corollary.
\qed

Let
\[L:=\{ x^\infty\mid \inf_{x\to x^\infty}\frac{Q(x)}{P(x)}>0\}.\]
Since \(\inf_{x\to x^\infty}\frac{Q(x)}{P(x)}>0\Leftrightarrow \sup_{x\to x^\infty}\frac{P(x)}{Q(x)}<\infty\), 
we see that there is a decreasing function \(h:\Nb\to \Qb\), i.e., \(\forall n< m\ h(n)\geq h(m)\),  such that 
\begin{equation}\label{assumpB}
P(L\cap\frac{P}{Q}>k)<h(k)\to 0\  (k\to\infty).
\end{equation}
We say that \(\frac{P}{Q}\) is  {\it effectively bounded in probability}  if there is a computable \(h\) in (\ref{assumpB}).

\begin{lemma}\label{lemma-equiv} Let \(P\) and \(Q\) be probabilities on \(\Omega\).\\
(a)  If \(\frac{Q}{P}\) is effectively-approximable and 
\(\frac{P}{Q}\) is effectively bounded in probability then
 \(\Hc^P\subseteq \Hc^{Q}\cup L^c\).\\ 
(b)  If \(P(L)=1\) and \(\frac{P}{Q}\) is effectively bounded  in probability, we have \(\Hc^P\subseteq \Hc^{Q}\). \\
(c) \(\exists c, c'\forall x\  0<c<\frac{Q(x)}{P(x)}<c'<\infty \Rightarrow \Hc^P=\Hc^Q.\)
\end{lemma}
Proof)
(a)
Let \(\{V_n\}\) be a test w.r.t.~\(Q\) and \(Q(\tilde{V}_n)<n^{-2}\) for all \(n\).
For \(n,m\in\Nb\), let 
\[L_{m,n}:=\{x\mid \frac{P(x)}{Q(x)}<m,\ |x|=n\}\text{ and }T_{m,n}:=\{x\mid f(x,n)<m,\ |x|=n\}.\]
From (\ref{r-approx}), there is a constant \(c\) and strictly increasing \(g\) such that 
\[\forall n,m\ T_{m-c,n}\subseteq L_{g^{-1}(m),n}\subseteq T_{m,n}.\]
Thus we have 
\[\cup_m\cap_n \tilde{L}_{m,n}=\cup_m\cap_n\tilde{T}_{m,n}=L.\]
Let \(V^m_n:=V_n\cap T_{m,n}\) then
\begin{align*}
P(\tilde{V}^m_n)&=  P(\tilde{V}^m_n\cap \{\frac{P}{Q}>k\})+P(\tilde{V}^m_n\cap \{\frac{P}{Q}\leq k\})\\
&\leq h(k)+kQ(\tilde{V}_n)< h(n)+n^{-1},
\end{align*}
where \(k=n\) in the last inequality. 
Since \(h\) is computable, there is a computable \(l:\Nb\to\Nb\) such that 
\(\sum_n h(l(n))+{l(n)}^{-1}<\infty\).
Since \(f\) is computable, \(\{V^m_{l(n)}\}_{n\in\Nb}\) is an uniformly r.e.~set and from Solovay's theorem (see \cite{shen89}),  we have
\[\limsup_n \tilde{V}^m_{l(n)}\subseteq (\Hc^P)^c.\]
Thus 
\[\cap_n\tilde{V}_n\cap (\cap_n\tilde{L}_{m,n})\subseteq \limsup_n \tilde{V}^m_n\subseteq \limsup_n \tilde{V}^m_{l(n)}\subseteq (\Hc^P)^c.\]
Since the above equation holds for all \(m\in\Nb\) and for any test \(V\) w.r.t.~\(Q\), we have
\[(\Hc^Q)^c\cap L=\cup_m\cup_{V: \text{ test w.r.t.}~Q} (\cap_n\tilde{V}_n\cap (\cap_n\tilde{L}_{m,n}))\subseteq (\Hc^P)^c.\]
Thus we have (a).\\
(b) 
Let \(\{V_n\}_{n\in\Nb}\) be an uniformly r.e.~set such that \(Q(\tilde{V}_n)<n^{-2}\) for all \(n\).
Since \(P(L)=1\), we have
\[P(\tilde{V}_n)\leq P(\tilde{V}_n\cap \frac{P}{Q}>k)+P(\tilde{V}_n\cap \frac{P}{Q}\leq k)< h(k)+kQ(\tilde{V}_n)=h(n)+n^{-1},\]
where \(k=n\).
Since \(h\) is computable, we see that  if \(\{V_n\}_{n\in\Nb}\) is a blind-test w.r.t.~\(Q\), it is a blind-test w.r.t.~\(P\).\\
(c) This follows immediately from (b). 
\qed

\begin{theorem}\label{th-classify}
Let \(P\) and \(Q\) be probabilities on \(\Omega\).
 If \(\frac{Q}{P}\) is effectively-approximable and 
\(\frac{P}{Q}\) is effectively bounded in probability then
\[
\Hc^P\cap\Hc^Q=\Hc^P\cap \{ x^\infty\mid 0<\inf_{x\to x^\infty}\frac{Q(x)}{P(x)}\},\text{ and }
\]
\[
\Hc^P\cap(\Hc^Q)^c=\Hc^P\cap \{ x^\infty\mid 0=\inf_{x\to x^\infty}\frac{Q(x)}{P(x)}\}.
\]
\end{theorem}
Proof)
Let
\(L:=\{ x^\infty\mid \inf_{x\to x^\infty}\frac{Q(x)}{P(x)}>0\}\).
Then
\begin{align}
\Hc^P\cap\Hc^Q & \subseteq \Hc^P\cap L\label{th-class-1}\\
&\subseteq (\Hc^P\cap L)\cap (\Hc^Q\cup L^c)\label{th-class-2}\\
&= \Hc^P\cap\Hc^Q\cap L\nonumber\\
&\subseteq \Hc^P\cap\Hc^Q\nonumber,
\end{align}
where (\ref{th-class-1}) follows from Corollary~\ref{col-q/p} and (\ref{th-class-2}) follows from Lemma~\ref{lemma-equiv} (a).
The second equation follows from the first one. 
\qed

Note that if \(P\) and \(Q\) are computable, we have 
\(\R^P\cap\R^{Q}=\R^P\cap\{x^\infty\mid 0<\inf_{x\sqsubset x^\infty} Q(x)/P(x)\}\), see \cite{{BM2007},{takahashiIandC2}}.
We can relative the results above, i.e., for any \(y^\infty\in\Omega\),  we can replace \(\Hc^P\) and \(\Hc^Q\) with \(\Hc^{P,y^\infty}\) and \(\Hc^{Q,y^\infty}\) in Corollary~\ref{lemma-equiv} and Theorem~\ref{th-classify}, respectively. 
If two conditions in Theorem~\ref{th-classify} are satisfied for conditional probabilities \(P(\cdot\vert y^\infty)\) and \(Q(\cdot\vert y^\infty)\),
we can replace the condition (vi) in Theorem~\ref{th-main} with \(x^\infty\in \Hc^{P(\cdot\vert y^\infty), y^\infty}\cap\R^{Q(\cdot\vert y^\infty), y^\infty}\).

\subsection{Example}
Next we show an example that holds equality in (\ref{generalB}) even if the conditional probability is not computable for all random parameters. 
\begin{theorem}\label{exampleA}
There is a computable probability \(P\) on \(X\times Y, X=Y=\Omega\) such that for all \(y^\infty\in\R^{P_Y}\),\\
(a) \(P(\cdot\vert y^\infty)\) is not computable  and (b)
 \(\R^P_{y^\infty}=\Hc^{P(\cdot\vert y^\infty),y^\infty}\).
\end{theorem}
Proof)
We construct a computable probability \(P\) on \(\Omega^2\) such that \(\R^P=\R^Q\), where \(Q\) is the product of uniform probabilities, i.e., \(Q(x,y)=Q_X(x)Q_Y(y)=2^{-(|x|+|y|)}\) for all \(x,y\in S\) and \(P_X=Q_X,\ P_Y=Q_Y\). 
Let \(e_1, e_2,\ldots\) be an enumeration of partial computable functions.
Let \(\Delta_1=\Delta(0), \Delta_2=\Delta(10), \Delta_3=\Delta(110),\ldots\) and we have \(\Omega\setminus \{1^\infty\}=\cup_n \Delta_n\).
We construct \(P\) such that for all \(n\) and \(y^\infty\in\R^{P_Y}\),
 \begin{equation}\label{oracle}
  P(\Delta_n| y^\infty)\ne P_X(\Delta_n)\Leftrightarrow e_n\text{ halts with oracle }y^\infty.
 \end{equation}
Observe that  there is a partial computable function \(U\) and a total computable function \(f\) such that for all \(n\) and \(y^\infty\), 
\begin{align}
  e_n\text{ halts with oracle }y^\infty\Leftrightarrow & \exists y\sqsubset y^\infty\ U(n,y)\text{ halts } \nonumber\\
\Leftrightarrow &  \exists y\sqsubset y^\infty\exists t\ f(n,y,t)=0 \nonumber\\
\Leftrightarrow  &\exists y\sqsubset y^\infty\ f(n,y, |y|)=0.\nonumber
\end{align}
Here if \(U(n,y)\) halts for some \(n,y\) then \(U(n,z)\) halts for all extension \(z\) of \(y\), and \(f(n,y,t)=0\) for some \(n,y,t\) then \(f(n,z,l)=0\) for all extension \(z\) of \(y\) and \(t\leq l\).
Intuitively,  the argument \(t\) of  \(f(n,y,t)\) is the number of steps  of the computation of  \(U(n,y)\). 
Let    \(P(x | \lambda)=P_X(x)=Q_X(x)=2^{-|x|}\) and \(P_Y(y)=Q_Y(y)=2^{-|y|}\) for all \(x, y\in S\).
Let \(0<\epsilon<1\), 
\[
P(\Delta_n|y0):=
\begin{cases}
P(\Delta_n|y)(1-\epsilon) &\text{ if }f(n,y, |y|)=0\text{ and }f(n,y', |y'|)\ne 0 \\
P(\Delta_n|y) &\text{ else, and}
\end{cases}
\]
\[
P(\Delta_n|y1):=
\begin{cases}
P(\Delta_n|y)(1+\epsilon) &\text{ if }f(n,y, |y|)=0\text{ and }f(n,y', |y'|)\ne 0 \\
P(\Delta_n|y) &\text{ else. }
\end{cases}
\]
Here \(|y|\) is the length of \(y\) and \(y'=y_1\ldots y_{|y|-1}\). 
By induction, we see that \(P(\Delta_n, y)=P(\Delta_n|y0)P(y0)+P(\Delta_n|y1)P(y1)\) for all \(n,y\).
If \(\Delta(x)\subseteq\Delta_n\) and \(|x|\geq n\) then let \(P(x| y)=2^{n-|x|}P(\Delta_n|y)\).
Then \(P\) is consistently defined, i.e., \(P(x,y)=\sum_{i,j\in\{0,1\}}P(xi,yj)\)  for all \(x,y\in S\), and \(P_X(x)=Q_X(x), P_Y(y)=Q_Y(y)\).
Since \(f\) is computable, we see that \(P\) is computable. 
By construction we have (\ref{oracle}) and the conditional probabilities are not computable for all \(y^\infty\in\R^{P_Y}\). 
We have 
\begin{align}
\forall x,y\ (1-\epsilon)\leq \frac{P(x,y)}{Q(x,y)}\leq (1+\epsilon), \label{equivA}\\ 
\forall x\ (1-\epsilon)\leq \frac{P(x|y^\infty)}{Q(x|y^\infty)}\leq (1+\epsilon).\label{equivB}
\end{align}
Hence
\begin{align}
(x^\infty, y^\infty)\in\R^P  & \Leftrightarrow (x^\infty, y^\infty)\in\R^Q \label{equivC}\\
& \Leftrightarrow y^\infty\in\R^{Q_Y}, x^\infty\in\R^{Q(\cdot\vert y^\infty),y^\infty} \label{equivD}\\
& \Leftrightarrow y^\infty\in\R^{P_Y}, x^\infty\in\Hc^{P(\cdot\vert y^\infty),y^\infty} \label{equivE},
\end{align}
where (\ref{equivC}) follows from Lemma~\ref{lemma-equiv} (c) and (\ref{equivA}), (\ref{equivD}) follows from Lambalgen's theorem, and (\ref{equivE}) follows from (relativized version of ) Lemma~\ref{lemma-equiv} (c) and (\ref{equivB}).
\qed

The following theorem shows an example  that equality does not hold in (\ref{generalB}). 
\begin{theorem}[Bauwens \cite{bauwens}]\label{exampleB}
There is a computable probability \(P\) on \(X\times Y, X=Y=\Omega\) and  \(y^\infty\in\R^{P_Y}\) such that\\
(a) \(P(\cdot\vert y^\infty)\) is not computable and 
(b)
 \(\R^P_{y^\infty}\ne\Hc^{P(\cdot\vert y^\infty),y^\infty}\).
\end{theorem}
\subsection{Other sufficient condition}
Finally we give another sufficient condition to the equality in  (\ref{generalB}), which requires a condition related to the convergence rate of conditional probability, 
however it does not require explicitly the existence of another computable conditional probability as in Theorem~\ref{th-main}.
\begin{lemma}\label{lem-A}
Let \(y^\infty\in \R^{P_Y}\).
If \(A\) is r.e.~relative to \(y^\infty\) and \(P(\tilde{A}\vert y^\infty)<\epsilon\) then
there are uniformly r.e.~sets \(U_1, U_2,\ldots\subseteq S\times S\) such that  
\begin{gather}
 \cup_n (\tilde{U}_n\smallsetminus \tilde{U}_{n+1})\cap \liminf_n \tilde{U}_n=\emptyset,\label{eq-decomp-sub}\\
 (\tilde{U}_n\smallsetminus \tilde{U}_{n+1})\cap (\tilde{U}_m\smallsetminus \tilde{U}_{m+1})=\emptyset\text{ for }n\ne m, \label{eq-decomp3}\\
 \widetilde{\cup_n U}  = \cup_n \tilde{U}_n  = \cup_n (\tilde{U}_n\smallsetminus \tilde{U}_{n+1})\cup \liminf_n \tilde{U}_n, \text{ and } \label{eq-decomp}\\
 \tilde{A}=(\liminf_n \tilde{U}_n)_{y^\infty}.\label{eq-decomp4}
\end{gather}
\end{lemma}
Proof)
Since \(A\) is r.e.~relative to \(y^\infty\), there is a partial computable~\(B: \Nb\times S\to S\) such that (i) if \(B(i,y)\) is defined then \(B(i,y)=B(i,y')\) for all \(y'\sqsupseteq y\) and (ii)
\(A=\{ x \mid \exists i\in\Nb, y\sqsubset y^\infty\ \ B(i,y)=x\}\).
Let \(T:=\{ (x,y)\mid \exists i\ \ B(i,y)=x\}\).
Since \(T\) is r.e., there is a non-overlapping r.e.~\(W\) such that \(\tilde{T}=\tilde{W}\).
Here \(W\) is called non-overlapping if \(\Delta(x,y)\cap\Delta(x',y')=\emptyset\) for \((x,y), (x',y')\in W\) and \((x,y)\ne (x',y')\), see \cite{takahashiIandC2}.
Let \(W:=\{(x^1,y^1), (x^2,y^2), \ldots\}\) be a recursive enumeration,  \(W_n:=\{(x^1,y^1),\ldots, (x^n,y^n)\}\), and 
\(V_n=\{y\in S\mid \Delta (y)\subseteq \cap_{1\leq i\leq n} C_i, \ C_i\in\{\Delta(y^i), \Delta(y^i)^c\}, 1\leq i\leq n\}\).
\(V_n\) is the partition generated from \(\Delta(y^i), i=1,\ldots,n\).
For \(A\subseteq S^2\), let \(A_y:=\{ x\mid \exists z\  (x,z)\in A, z\sqsubseteq y\}\).
For example, \(W_{n,y}=\{x\mid \exists z\  (x,z)\in W_n, z\sqsubseteq y\}\).
Set
\begin{equation}\label{def-U}
U_n:=\{ (x,y)\mid \sum_{x\in W_{n, y}}P(x | y)<\epsilon,\ y\in V_n\}.  
\end{equation}
Let \(x\in A\smallsetminus B\leftrightarrow x\in A  \land x\notin B\) for sets \(A\) and \(B\).
Since \(\liminf_n \tilde{U}_n= \cup_m \cap_{n\geq m}\tilde{U}_n\), we have (\ref{eq-decomp-sub}) and (\ref{eq-decomp}).

Next we show that for \(n,m\in\Nb\),
\begin{gather}
\exists \text{ open set }O_n\  \ \tilde{U}_n = \tilde{U}_n\cap (\Omega\times O_n)\label{decomp2-1},\\
\tilde{U}_n\smallsetminus \tilde{U}_{n+1}=\tilde{U}_n\cap (\Omega\times D_n),\text{ where }D_n:=O_n\smallsetminus O_{n+1},\label{decomp2-2}\text{ and }\\
D_n\cap D_m=\emptyset\text{ for }n\ne m\label{decomp2-3}.
\end{gather}
Let  \(O_n:=\{y^\infty\mid \tilde{U}_{n,y^\infty}\ne\emptyset\}\) for \(n\in\Nb\) then from (\ref{def-U}) and \(P(\cdot\vert y)\) converges as \(y\to y^\infty\in\R^{P_Y}\) (see \cite{takahashiIandC}), we see that \(O_n\) is open.
Since \(W_n\subseteq W_m\) for \(n<m\), we have 
\begin{equation}\label{sub-equiv}
y^\infty\in O_n\cap O_m\Leftrightarrow \tilde{U}_{m,y^\infty}\supseteq\tilde{U}_{n,y^\infty}\ne\emptyset,
\end{equation}
and
(\ref{decomp2-1}) and (\ref{decomp2-2}) hold. 
Suppose that \(D_n\cap D_m\ne\emptyset\) for \(n<m\).
Since \(y^\infty\in D_n\cap D_m\subseteq O_n\cap O_m\), we have
\(\emptyset\ne W_{n,y}\subseteq W_{n+1,y}\subseteq W_{m, y}\).
Thus we have \(y^\infty\in O_{n+1}\) and  \(y^\infty\notin D_n\), which is a contradiction, and  we have (\ref{decomp2-3}) and (\ref{eq-decomp3}).

Since \(\lim_{y\to y^\infty}P(\tilde{B} | y)=P(\tilde{B} | y^\infty)\) for finite set \(B\) and   \(y^\infty\in\R^{P_Y}\), we have if \(P(\tilde{B} | y^\infty)<\epsilon\) then there is \(y\sqsubset y^\infty\)
such that \(P(\tilde{B} | y)<\epsilon\).
Thus we have \(A\subseteq (\liminf_n \tilde{U}_n)_{y^\infty}\subseteq (\cup_n \tilde{U}_n)_{y^\infty}\subseteq \tilde{W}_{y^\infty}=A\), and we have (\ref{eq-decomp4}).
\qed

From (\ref{eq-decomp3}), we see that there  is \(f:\Qb\to \Nb\) such that 
\begin{equation}\label{lem-condition}
\forall \epsilon>0\ P(\cup_{f(\epsilon)\leq n} \tilde{U}_n\smallsetminus \tilde{U}_{n+1})<\epsilon.
\end{equation}

\begin{theorem}
Fix \(y^\infty\in \R^{P_Y}\). 
Assume that  there is a computable \(f\) that satisfies (\ref{lem-condition}) for any \(\epsilon'>0\) and r.e.~relative to \(y^\infty\), \(A\) such that \(P(\tilde{A}\vert y^\infty)<\epsilon'\) in Lemma~\ref{lem-A}.
Then
\[\R^P_{y^\infty}=\Hc^{P(\cdot\vert y^\infty),y^\infty}.\]
\end{theorem}
Proof)
If \(A\) is r.e.~relative to \(y^\infty\) such that \(P(\tilde{A}\vert y^\infty)<\epsilon'\) then 
\(\cup_{f(\epsilon)\leq n} U_n\) is r.e.~and
\(P(\cup_{f(\epsilon)\leq n} \tilde{U}_n)=P(\cup_{f(\epsilon)\leq n} (\tilde{U}_n\setminus \tilde{U}_{n+1}) \cup \liminf_n \tilde{U}_n)<2\epsilon'\).
Thus we see that if \(A\) is a test w.r.t.~\(P(\cdot\vert y^\infty)\) then it is covered by a test w.r.t.~\(P\).
\qed

\begin{center}
Acknowledgement
\end{center}
This paper is based on the work when the author visited LIRMM Montpellier France. 
The author thanks  Prof.~A.~Shen (LIRMM France),  Prof.~A.~Romashchenko (LIRMM France), Prof.~Bruno Bauwens (Nancy France), Prof.~Teturo Kamae (Osaka city univ.), Prof.~Hiroshi Sugita (Osaka univ.), and Prof.~Akio Fujiwara (Osaka univ.) for discussions and comments. 
A part of work is supported by JSPS KAKENHI Grant number 24540153.

\end{document}